\documentstyle{amsppt}
\magnification=\magstep1
\hsize=30truecc
\baselineskip=16truept

\topmatter

\title
Partial Measures
\endtitle

\author
O.~E.~Tikhonov
\endauthor

\address
Research Institute of Mathematics and Mechanics, 
Kazan State University, Universitetskaya Str.\ 17, Kazan, Tatarstan, 
420008 Russia
\endaddress

\email
Oleg.Tikhonov{\@}ksu.ru
\endemail

\keywords
Set function, $\sigma$-algebra, Jordan decomposition theorem, 
Radon-Nikodym theorem
\endkeywords 

\subjclass
28A10
\endsubjclass

\thanks
Supported by the Russian Foundation for Basic 
Research (grant no.\ 98--01--00103) and the scientific program 
Universities of Russia --- Basic Research (grant no.\ 1617).
\endthanks

\abstract
We study $\sigma$-additive set functions defined on a hereditary 
subclass of a $\sigma$-algebra and taken values in the extended real 
line. Analogs of the Jordan decomposition theorem and the 
Radon-Nikodym theorem are obtained.
\endabstract

\endtopmatter

Throughout this note, ${\Cal A}$ stands for a $\sigma$-algebra of 
subsets of a set $\Omega$, $\overline{\Bbb R}$ denotes the extended 
real line ${\Bbb R} \cup \{ -\infty \} \cup \{ +\infty \}$. 
By a {\it measure} on ${\Cal A}$ we mean a $\sigma$-additive set 
function $\mu : {\Cal A} \to \overline{\Bbb R}$ 
with $\mu (\emptyset ) = 0$. Note that in fact a measure can take at 
most one of the values $+\infty$, $-\infty$.

\definition {Definition} 
A set function $\mu : {\Cal D } (\mu ) \to \overline{\Bbb R}$ defined 
on a nonvoid subclass ${\Cal D } (\mu )$ of a $\sigma$-algebra 
$\Cal A$ is called a {\it partial measure} on ${\Cal A}$ if for every 
$B \in {\Cal D } (\mu )$ 
the restriction $\mu |_{{\Cal A} \cap B}$ is a measure on the 
$\sigma$-algebra ${\Cal A} \cap B = \{ A \cap B : A \in {\Cal A} \}$.
A partial measure is called {\it maximal} if it has no proper 
extension to a partial measure on ${\Cal A}$. 
\enddefinition

By making use of the Zorn lemma, it is easy to see that every partial 
measure can be extended to a maximal one. 

\remark{Example 1}
Let $\mu _1$, $\mu _2$ be positive measures on a ${\Cal A}$. Then the 
formula $(\mu _1 - \mu _2) (A) = \mu _1 (A) - \mu _2 (A)$ 
defines a partial measure $\mu _1 - \mu _2$. 
(${\Cal D } (\mu _1 - \mu _2) = 
\{ A \in {\Cal A} : \mu _1 (A) - \mu _2 (A)
\text{ is well-posed} \}$.) 
\endremark

\remark{Example 2} 
Let $(\Omega , {\Cal A} , {\Bbb P})$ be a probability space, 
$\xi : \Omega \to \overline{\Bbb R}$ be a random variable. Then the 
formula $\mu _{\xi} (A) = \int _A \xi \, d {\Bbb P}$ defines a maximal 
partial measure $\mu _{\xi}$ on ${\Cal A}$. 
(${\Cal D } (\mu _{\xi} ) = \{ A \in {\Cal A} : \xi \text{ is 
quasiintegrable over } A \}$.) 
\endremark

For a partial measure $\mu$ on $\Cal A$, define two classes: 
${\Cal F}^+ (\mu ) = \{ F \in {\Cal D} (\mu ) : 
{\Cal A} \ni A \subset F \implies \mu (A) \ge 0 \}$,  
${\Cal F}^- (\mu ) = \{ F \in {\Cal D} (\mu ) : 
{\Cal A} \ni A \subset F \implies \mu (A) \le 0 \}$.

\proclaim{Lemma 1}
Let $\{ F_i : i \in I \}$ and $\{ F'_j :j \in J \}$ be two 
finite or countable families of mutually disjoint sets of 
${\Cal F}^+ (\mu )$ such that 
$\bigcup _{i \in I} F_i = \bigcup _{j \in J} F'_j $. Then 
$\sum _{i \in I} \mu (F_i) = \sum _{j \in J} \mu (F'_j)$. 
\endproclaim

\demo{Proof}
Since 
$F_i = \bigcup _{j \in J} (F_i \cap F'_j)$ for every $i \in I$,
$F'_j = \bigcup _{i \in I} (F_i \cap F'_j)$ for every $j \in J$, 
and $F_i , F'_j \in {\Cal D} (\mu )$ we have
$$
\sum _{i \in I} \mu (F_i) 
= \sum _{i \in I} \sum _{j \in J} \mu (F_i \cap F'_j)
= \sum _{j \in J} \sum _{i \in I} \mu (F_i \cap F'_j)
= \sum _{j \in J} \mu (F'_j) .
$$
\enddemo
                                
\proclaim{Lemma 2}
Let a partial measure $\mu$ on $\Cal A$ be maximal. 
If $\{ F_i : i \in I \}$ is a finite or countable family of sets of 
${\Cal F}^+ (\mu )$ then 
$\bigcup _{i \in I} F_i \in {\Cal F}^+ (\mu )$.
\endproclaim

\demo{Proof}
Clearly, it suffices to prove that all the unions 
$\bigcup _{i \in I} F_i$ of finite or countable 
families of mutually disjoint sets of ${\Cal F}^+ (\mu )$ belong to 
${\Cal F}^+ (\mu )$. Extend $\mu $ to such unions by means of the 
formula 
$
\widetilde{\mu} \bigl ( \bigcup _{i \in I} F_i \bigr ) = 
\sum _{i \in I} \mu ( F_i ) 
$.
(By Lemma 1, $\widetilde{\mu}$ is well-defined.) The extension obtained 
is a partial measure on ${\Cal A}$. Really, if 
$A = \bigcup _{i \in I} F_i$ ($I$ is finite or countable, 
$F_i \in {\Cal F}^+ (\mu )$ for any $i \in I$, 
$F_i \cap F_j = \emptyset$ for $i \ne j$) and  
$A = \bigcup _{k \in K} A_k$ ($K$ is finite or countable, 
$A_k \in {\Cal A}$ for any $k \in K$, 
$A_k \cap A_l = \emptyset$ for $k \ne l$) then 
$$
\align
\widetilde{\mu} (A) & = 
\sum _{i \in I} \mu \bigl ( \bigcup _{k \in K} (F_i \cap A_k ) \bigr ) 
= \sum _{i \in I} \sum _{k \in K} \mu (F_i \cap A_k ) = 
\sum _{k \in K} \sum _{i \in I} \mu (F_i \cap A_k ) \\
& = \sum _{k \in K} \widetilde{\mu} (A_k) .
\endalign
$$
Moreover, it is easy to check that the unions under consideration are 
sets of the class ${\Cal F}^+ (\widetilde{\mu})$. To complete the 
proof, it suffices to observe that $\widetilde{\mu} = \mu$ by the 
maximality of $\mu$.
\enddemo

\proclaim{Theorem 1}
Let $\mu$ be a maximal partial measure on a $\sigma$-algebra $\Cal A$. 
Then the formulas
$$
\mu ^+ (A) = \sup _{F \in {\Cal F}^+ (\mu )} \mu (A \cap F) , \quad 
\mu ^- (A) = \sup _{F \in {\Cal F}^- (\mu )} (- \mu (A \cap F)) \quad 
(A \in {\Cal A}) 
$$
define two positive measures on $\Cal A$ such that 
$\mu = \mu ^+ - \mu ^-$. 
Moreover, if $\mu _1$ is a positive measure on $\Cal A$ with 
$\mu (A) \le \mu _1 (A)$ for all $A \in {\Cal D} (\mu )$ then  
$\mu ^+ \le \mu _1$. 
Similarly, if $\mu _2$ is a positive measure on $\Cal A$ with 
$- \mu (A) \le \mu _2 (A)$ for all $A \in {\Cal D} (\mu )$ then  
$\mu ^- \le \mu _2$. 
\endproclaim

\demo{Proof}
Observe, first, that for every $F \in {\Cal F}^+ (\mu )$ the formula 
$\mu _F (A) = \mu (A \cap F)$ defines a positive measure $\mu _F$ on 
$\Cal A$. 
Also, from Lemma 2 it follows that the class ${\Cal F}^+ (\mu )$ is 
upwards directed by inclusion. 
Moreover, the condition $F_1 \subset F_2$ 
($F_1 , F_2 \in {\Cal F}^+ (\mu )$) implies  
$\mu _{F_1} (A) \le \mu _{F_2} (A)$ for all $A \in {\Cal A}$. 
It follows that the set function $\mu ^+$ on $\Cal A$ is a pointwise 
limit of an increasing net of positive measures, therefore, $\mu ^+$ 
is a positive measure in turn. 
Similarly, $\mu ^-$ is a positive measure on $\Cal A$.

Next, let $A \in {\Cal D} (\mu )$. By the Hahn decomposition theorem, 
$A$ can be represented as a union, $A = A^+ \cup A^-$, of two disjoint sets 
$A^+ \in {\Cal F}^+ (\mu )$ and $A^- \in {\Cal F}^- (\mu )$. Clearly, 
$\mu ^+ (A) = \mu (A^+ )$ and $\mu ^- (A) = - \mu (A^-)$. Thus, 
$\mu (A) = \mu ^+ (A) - \mu ^- (A)$ for any $A \in {\Cal D} (\mu )$. 
Since $\mu ^+ - \mu ^-$ is a partial measure on $\Cal A$ (see Example 
1) the maximality of $\mu$ yields $\mu = \mu ^+ - \mu ^-$.

Now, let $\mu _1$ is a positive measure on $\Cal A$ with 
$\mu (A) \le \mu _1 (A)$ for all $A \in {\Cal D} (\mu )$. 
Take $A \in {\Cal A}$. For every $F \in {\Cal F}^+ (\mu )$ we have
$\mu (A \cap F) \le \mu _1 (A \cap F) \le \mu _1 (A)$. Therefore, 
$\mu ^+ (A) = 
\sup \limits _{F \in {\Cal F}^+ (\mu )} \mu (A \cap F) 
\le \mu _1 (A)$. 
Thus, we have proved the extremal property for $\mu ^+$. 
The similar extremal property for $\mu ^-$ can be proved in the same 
manner.
\enddemo

\proclaim{Corollary 1}
Let $\mu$ be a maximal partial measure on a $\sigma$-algebra $\Cal A$ 
and $A$ be a set of ${\Cal A} \setminus {\Cal D} (\mu )$. 
Then there exist two subsets $A' , A''$ of $A$ such that 
$A' \in {\Cal F}^+ (\mu )$ and $\mu (A') = + \infty$, 
$A'' \in {\Cal F}^- (\mu )$ and $\mu (A'') = - \infty$.
\endproclaim

\demo{Proof}
Since $A \not \in {\Cal D} (\mu ) = {\Cal D} (\mu ^+ - \mu ^-)$ we 
have $\mu ^+ (A) = + \infty$ and $\mu ^- (A) = + \infty$. By the 
definition of $\mu ^+$, we can find a sequence 
$\{ A_n \} \subset {\Cal F}^+ (\mu )$ of subsets of $A$, 
such that $\mu (A_n) \ge n$. 
By Lemma 2, the union $A' = \bigcup _{n=1} ^{\infty} A_n$ is a set of  
${\Cal F}^+ (\mu )$, and, obviously, $\mu (A') = + \infty$. 
The set $A''$ can be constructed in the same manner.
\enddemo

\remark{Remark}
The preceding theorem is an analog of the Jordan decomposition 
theorem. Example 3 below shows that the direct analog of the Hahn 
decomposition does not hold for general maximal partial measures. 
Nevertheless, the latter analog does hold provided that a maximal 
partial measure under consideration is absolutely continuous with 
respect to a certain ``good'' measure (see the proof of Theorem 
2 below). 
\endremark

\remark{Example 3}
Let a $\sigma$-algebra $\Cal A$ of subsets of $\Omega$ contains all 
the one-point subsets of $\Omega$ and does not coincide with the class 
${\Cal P} (\Omega )$ of all the subsets of $\Omega$. (For instance, we 
may take the Borel $\sigma$-algebra of $[0,1]$ as $\Cal A$.) 
Let $B \in {\Cal P} (\Omega ) \setminus {\Cal A}$. 
Define the set function $\mu $ as follows:
$$
\mu (A) = \cases 
0 & \text{if $A = \emptyset$ ,} \\
+ \infty & \text{if ${\Cal A} \ni A \subset B$ and $A \ne \emptyset$ ,} \\ 
- \infty & 
\text{if ${\Cal A} \ni A \subset \Omega \setminus B$ and $A \ne \emptyset$ .} 
\endcases
$$
It is easy to see that $\mu $ is a maximal partial measure on $\Cal A$ 
and that there exists no $C \in {\Cal F} ^+ ( \mu )$ with 
$\Omega \setminus C \in {\Cal F} ^- ( \mu )$.
\endremark

In the rest of the note, the triple $(\Omega , {\Cal A} , {\Bbb P})$ 
stands for a probability space, the abbreviation a.\,s.\ means ``almost 
surely with respect to ${\Bbb P}$''. 
The following theorem is an analog of the Radon-Nikodym one.

\proclaim{Theorem 2} 
Let a maximal partial measure $\mu$ on $\Cal A$ be absolutely 
continuous with respect to a probability measure ${\Bbb P}$ 
(i.\,e., ${\Bbb P} (A) = 0$ implies $A \in {\Cal D} ( \mu )$ and 
$\mu (A) = 0$). Then there exists an a.\,s.\ unique random variable 
$\xi : \Omega \to \overline{\Bbb R}$ such that $\mu = \mu _{\xi}$, 
where $\mu _{\xi} (A) = \int _A \xi \, d {\Bbb P}$ (see Example 2).
\endproclaim

\demo{Proof}
Set $\Omega ^+ = \operatorname{ess} \, \sup {\Cal F}^+ ( \mu )$. 
Recall that $\operatorname{ess} \, \sup {\Cal F}^+ ( \mu )$ is a set 
of $\Cal A$ such that for any $A \in {\Cal A}$
$$
\forall F \in {\Cal F}^+ ( \mu ) \,
(F \underset \text{a.\,s.} \to \subset A) \iff 
\operatorname{ess} \, \sup {\Cal F}^+ ( \mu ) 
\underset \text{a.\,s.} \to \subset A 
$$
(see, e.\,g., [1, Proposition II.4.1]). From the proof of the 
proposition just referred, it follows that $\Omega ^+$ can be 
constructed as a union 
%$\bigcup _{i \in I} F_i$ 
of a countable family of sets of ${\Cal F}^+ ( \mu )$. 
Therefore, $\Omega ^+ \in {\Cal F}^+ ( \mu )$ by Lemma 2.  

Our next goal is to show that $\Omega ^- = \Omega \setminus \Omega ^+$ 
is a set of ${\Cal F}^- ( \mu )$. Suppose, on the contrary, that  
$\Omega ^- \not \in {\Cal F}^- ( \mu )$. Then there exists 
$A \subset \Omega ^-$ ($A \in {\Cal A})$ such that either $\mu (A) > 0$ 
or $A \not \in {\Cal D} ( \mu )$. In each of the cases, we can find a 
subset $A' \subset A$ with $A' \in {\Cal F}^+ ( \mu )$ and 
$\mu (A') > 0$ (see Corollary 1 for the second case). 
Hence ${\Bbb P} (A') > 0$. 
This contradicts the definition of 
$\Omega ^+ = \operatorname{ess} \, \sup {\Cal F}^+ ( \mu )$.

By making use of the Radon-Nikodym theorem (see, e.\,g., 
[1, Proposition IV.1.4]), we can find two positive random variables 
$\xi ^+$ and $\xi ^-$ such that 
$\xi ^+ ( \omega ) = 0$ for $\omega \in \Omega ^-$, 
$\xi ^- ( \omega ) = 0$ for $\omega \in \Omega ^+$, 
$\mu (A) = \int _A \xi ^+ \, d {\Bbb P}$ if 
$\Omega ^+ \supset A \in {\Cal A}$, 
$\mu (A) = - \int _A \xi ^- \, d {\Bbb P}$ if 
$\Omega ^- \supset A \in {\Cal A}$. 
Writing $\xi = \xi ^+ - \xi ^-$ we have for $A \in {\Cal D} ( \mu )$
$$
\mu (A) = \mu (A \cap \Omega ^+ ) + \mu (A \cap \Omega ^- ) 
= \int _{A \cap \Omega ^+} \xi ^+ \, d {\Bbb P} - 
\int _{A \cap \Omega ^-} \xi ^- \, d {\Bbb P} = 
\int _A \xi \, d {\Bbb P} . 
$$
Hence, $\mu = \mu _{\xi}$ by the maximality of $\mu$.

It remains to prove the claim about the uniqueness. Let $\eta$ be a 
random value such that $\mu _{\eta} = \mu$. Then, by the uniqueness 
assertion in the Radon-Nikodym theorem, we have 
$\eta |_{\Omega ^+} \underset \text{a.\,s.} \to = \xi |_{\Omega ^+}$ 
and 
$\eta |_{\Omega ^-} \underset \text{a.\,s.} \to = \xi |_{\Omega ^-}$, 
hence $\eta \underset \text{a.\,s.} \to = \xi$.   
\enddemo


\Refs

\ref \no 1  
\by J. Neveu 
\book Bases math{\'e}matiques du calcul des probabilit{\'e}s
\publ Masson et Cie
\publaddr Paris
\yr 1964 
\transl Russian transl.
\publ Mir 
\publaddr Moscow
\yr 1969
\endref

\endRefs
\enddocument